\documentclass{gtart_h}

\def\ifplaintex{\expandafter\ifx\csname documentclass\endcsname\relax}

\ifplaintex 
\else
\fi

\def\gt{{\mathsurround=0pt\it $\cal G\mskip-2mu$eometry \&\ 
$\cal T\!\!$opology}}        

\def\gtm{{\mathsurround=0pt\it $\cal G\mskip-2mu$eometry \&\ 
$\cal T\!\!$opology $\cal M\mskip-1mu$onographs}}    

\def\gtp{{\mathsurround=0pt\it $\cal G\mskip-2mu$eometry \&\ 
$\cal T\!\!$opology $\cal P\!$ublications}}  

\def\recd{{\small Received:\qua\receiveddate\ifx\reviseddate\relax
\else\qquad Revised:\qua\reviseddate\fi\par}} 

\def\volumenumber#1{\def\thevolumenumber{#1}}
\def\volumeyear#1{\def\thevolumeyear{#1}}
\def\volumename#1{\def\thevolumename{#1}}
\def\papernumber#1{\def\thepapernumber{#1}}
\def\pagenumbers#1#2{\def\startpage{#1}\def\finishpage{#2}}
\def\published#1{\def\publishdate{#1}}
\def\received#1{\def\receiveddate{#1}}
\def\revised#1{\def\reviseddate{#1}}
\def\accepted#1{\def\accepteddate{#1}}
\def\asciititle#1{\def\theasciititle{#1}}

\def\asciiaddress#1{\def\theasciiaddress{#1}}
\def\asciiemail#1{\def\theasciiemail{#1}}

\long\def\asciiabstract#1{\long\def\theasciiabstract{#1}}
\def\asciikeywords#1{\def\theasciikeywords{#1}}

\let\\\par
\let\thevolumenumber\relax\let\thepapernumber\relax
\let\thevolumeyear\relax\let\startpage\relax
\let\finishpage\relax\let\publishdate\relax\let\receiveddate\relax
\let\reviseddate\relax\let\accepteddate\relax\let\theasciititle\relax
\let\theasciiauthors\relax\let\theasciiaddress\relax
\let\theasciiabstract\relax\let\theasciikeywords\relax

\let\theerratum\relax\let\theasciiemail\relax
\let\theshortauthors\relax\let\theshorttitle\relax

\def\startpage{1}\def\finishpage{15}\def\thepapernumber{77}

\volumenumber{2}
\volumename{Proceedings of the Kirbyfest}
\volumeyear{1999}

\long\def\maketitlep{   

\count0=\startpage

\gtm\nl        
{\small Volume \thevolumenumber: \thevolumename\nl 
\ifx\theerratum\relax\else Erratum \erratumnumber\nl\fi
Pages \startpage--\finishpage\nl}

\vglue 0.1truein   

{\parskip=0pt\leftskip 0pt plus 1fil\def\\{\par\smallskip}{\ifplaintex\large
\else\Large\fi\bf\thetitle}\par\medskip}   
\vglue 0.05truein 

{\parskip=0pt\leftskip 0pt plus 1fil\def\\{\par}{\sc\theauthors}
\par\medskip}%
 
\vglue 0.03truein 

{\small\leftskip 25pt\rightskip 25pt{\bf Abstract}\stdspace\theabstract

{\bf AMS Classification}\stdspace\theprimaryclass
\ifx\thesecondaryclass\relax\else; \thesecondaryclass\fi\par
{\bf Keywords}\stdspace \thekeywords\par}\vglue 7pt

}   

\font\phead=cmsl9 scaled 950
\font=cmsl9 scaled 1050
\font\pnum=cmbx10 scaled 913
\font\lnum=cmbx10 
\font\pfoot=cmsl9 scaled 950
\font=cmsl9 scaled 1050
\ifplaintex
\headline{\vbox to 0pt{\vskip -4.5mm\line{\small\phead\ifnum
\count0=\startpage ISSN 1464-8997 (on line)
1464-8989 (printed) \hfill {\pnum\folio}\else\ifodd\count0\def\\{ }%
\ifx\theshorttitle\relax\thetitle\else\theshorttitle\fi\hfill{\pnum\folio}
\else\def\\{ and }{\pnum\folio}\hfill\ifx\theshortauthors\relax\theauthors
\else\theshortauthors\fi\fi\fi}\vss}}
\footline{\vbox to 0pt{\vglue 0mm\line{\small\pfoot\ifnum\count0=\startpage
Published \publishdate:\qua\copyright\ \gtp\hfill\else
\gtm, Volume \thevolumenumber\ (\thevolumeyear)\hfill\fi}\vss
}}
\else
\makeatletter
\def\@oddhead{{\small\lhead\ifnum\count0=\startpage ISSN 1464-8997 (on line)
1464-8989 (printed) \hfill {\lnum\number\count0}\else\ifodd\count0
\def\\{ }\ifx\theshorttitle\relax \thetitle \else\theshorttitle\fi\hfill
{\lnum\number\count0}\else\def\\{ and }{\lnum\number\count0}
\hfill\ifx\theshortauthors\relax 
\theauthors\else\theshortauthors\fi\fi\fi}}\def\@evenhead{@oddhead}
\def\@oddfoot{\small\lfoot\ifnum\count0=\startpage Published \publishdate:\qua\copyright\ \gtp\hfill\else
\gtm, Volume \thevolumenumber\ (\thevolumeyear)\hfill\fi}
\def\@evenfoot{@oddfoot}
\makeatother
\fi

\let\maketitlepage\maketitlep
\let\makeshorttitle\maketitlepage
\let\maketitle\maketitlepage

\newwrite\gtoutfile
\long\gdef\makeheadfile{  
{\def\\{, }\def\s{ }
\immediate\openout\gtoutfile head.xxx
\immediate\write\gtoutfile{Proxy-for: \ifx\theasciiauthors\relax
\theauthors\else\theasciiauthors\fi\s<\ifx\theasciiemail\relax\theemail\else\theasciiemail\fi>}
\immediate\write\gtoutfile{\noexpand\\}
\immediate\write\gtoutfile{Authors: \ifx\theasciiauthors\relax
\theauthors\else\theasciiauthors\fi}
{\def\\{ }\immediate\write\gtoutfile{Title: \ifx\theasciititle\relax
\thetitle\else\theasciititle\fi}}
\immediate\write\gtoutfile{Subj-class: GT or SG, GR etc}
\immediate\write\gtoutfile{MSC-class: \theprimaryclass\ifx\thesecondaryclass\relax\else, \thesecondaryclass\fi}
\immediate\write\gtoutfile{Journal-ref: Geom. Topol. Monogr. \thevolumenumber\s
(\thevolumeyear) \startpage-\finishpage}
\immediate\write\gtoutfile{Comments: Published by Geometry and Topology Monographs at}
\immediate\write\gtoutfile{\s\s\s  http://www.maths.warwick.ac.uk/gt/GTMon\thevolumenumber/paper\thepapernumber.abs.html}
\immediate\write\gtoutfile{\noexpand\\}
\immediate\write\gtoutfile{}
\ifx\theasciiabstract\relax
\immediate\write\gtoutfile{\theabstract}\else
\immediate\write\gtoutfile{\theasciiabstract}\fi
\immediate\write\gtoutfile{}
\immediate\write\gtoutfile{\noexpand\\}
\immediate\write\gtoutfile{}
\immediate\closeout\gtoutfile}}  

\def\maketitlepage{\maketitlep\makeheadfile}
\let\makeshorttitle\maketitlepage
\let\maketitle\maketitlepage

\volumenumber{7}
\volumename{Proceedings of the Casson Fest}
\volumeyear{2004}
\papernumber{13}
\pagenumbers{311}{333}
\received{31 October 2003}
\revised{21 March 2004}
\accepted{2 April 2003}
\published{20 September 2004}

\usepackage{amsmath,amssymb,amscd,curves}

\theoremstyle{plain}
\newtheorem{thm}{Theorem}[section]
\newtheorem{lem}[thm]{Lemma}

\newtheorem{prop}[thm]{Proposition}

\newtheorem*{ther}{Theorem}

\theoremstyle{definition}

\def \CPb {\overline{\mathbb{CP}}^{\,2}}
\def \CP {{\mathbb{CP}}^{\,2}} 
\def \R {\mathbb{R}}
\def \Z {\mathbb{Z}}
\def \Sig{\Sigma}

\def \vt {\vartheta}
\def \vp {\varphi}

\def \a {\alpha}
\def \b {\beta}
\def \g {\gamma}
\def \d {\delta}

\def \lam {\lambda}
\def \G {\Gamma}

\def \o {\omega}
\def \s {\sigma}
\def \t {\tau}
\def \z {\zeta}

\def \bg {\bar{\gamma}}
\def \bd {\partial}
\def \x {\times}
\def \- {\setminus}
\def \C {\subset}

\def\pt{\text{pt}}
\newcommand{\spin}{\operatorname{Spin}}

\def \ssw {\text{SW}}
\def \sw {\mathcal{SW}}
\def \swb{\overline {\mathcal{SW}}}

\def \DD {\Delta}

\def \Ds {\DD^{\text{\it{sym}}}}

\def \DN {\nabla}
\def\DNS{\nabla^{\text{\it{sym}}}}

\def \tz {\widetilde{\z}}
\def \X {\widetilde{X}}
\def \f {\bar{f}}
\begin{document}

\title{Tori in symplectic 4--manifolds}
\asciititle{Tori in symplectic 4-manifolds}

\authors{Ronald Fintushel\\Ronald J Stern}

\address{Department of Mathematics, Michigan State 
University\\East Lansing, Michigan 48824, USA}
\secondaddress{Department of Mathematics, University of
California\\Irvine, California 92697, USA}

\asciiaddress{Department of Mathematics, Michigan State
University\\East Lansing, Michigan 48824, USA\\and\\Department of 
Mathematics, University of California\\Irvine, California 92697, USA}

\gtemail{\mailto{ronfint@math.msu.edu}, \mailto{rstern@math.uci.edu}}
\asciiemail{ronfint@math.msu.edu, rstern@math.uci.edu}

\begin{abstract}
We study the question of how many embedded symplectic or Lagrangian tori
can represent the same homology class in a simply connected symplectic
4--manifold.
\end{abstract}
\asciiabstract{%
We study the question of how many embedded symplectic or Lagrangian
tori can represent the same homology class in a simply connected
symplectic 4-manifold.}


\asciikeywords{4-manifold, Seiberg-Witten invariant, symplectic, Lagrangian}

\maketitle

\section{Introduction}

The basic question addressed in this paper is: 

\begin{quotation}
\noindent Let X be a simply connected symplectic 4--manifold and let
$x \in H_2(X, \Z)$. How unique is an embedded symplectic or Lagrangian
representative of $x$?
\end{quotation}

It is only in the last few years that an answer to this question
has begun to emerge. The answer is `not very' for symplectic  tori of
self-intersection $0$ and remains elusive for higher genus surfaces. As
we show below:

\begin{ther}
If $X$ is a simply connected symplectic 4--manifold containing an
embedded symplectic torus $T$ of self-intersection $0$, then  for each
fixed integer $m\ge2$, there are infinitely many embedded symplectic
tori, each representing the homology class $m[T]$, no two of which
are equivalent under a smooth isotopy  of $X$.
\end{ther}   

The first such examples were produced by the present authors in
\cite{surfaces}, and the technique therein was enhanced to produce
further examples in \cite{V1, EP}. In section \ref{FibS} we give a proof of
the above theorem. In section \ref{sc} we give a proof of the theorem for
even $m\ge 6$ which is straightfoward, and which depends on some nice
theorems of Montesinos and Morton \cite{MM} and of Kanenobu \cite{K}
rather than on explicit constructions.

Some intriguing questions remain. Siebert and Tian have conjectured that
for symplectic 4--manifolds with $b^+=1$ and $c_1^2>0$ any embedded
symplectic surface must be symplectically isotopic to a holomorphic
curve. (Of course, no such manifold contains an embedded symplectic
torus of square $0$.) They have shown that in $\CP$ this is true for
each curve of degree $\le17$, and they have also some results to this
effect in $S^2\x S^2$. However the general problem is still wide open,
as is the case in general for surfaces of higher genus or for other
self-intersections. However, in the case where $\pi_1(X)\ne 0$ Ivan
Smith has constructed examples of nonisotopic but homologous surfaces
of square $0$ distinguished by $\pi_1$ of their complements.

Much less is known in the case of Lagrangian tori. Until this year,
it was unknown if there existed Lagrangian tori which were homologous
but inequivalent (under either isotopy or orientation-preserving
diffeomorphism). The first examples are due to Stefano Vidussi:

\begin{ther}{\rm\cite{V}}\qua
Let $K$ denote the trefoil knot. Then in the symplectic manifold $E(2)_K$
there is a primitive homology class $\a$ so that for each positive integer
$m$, there are infinitely many embedded Lagrangian tori representing
$m\a$, no two of which are equivalent under orientation-preserving
diffeomorphisms.
\end{ther}

Utilizing an invariant coming from Seiberg--Witten theory and the geometry
of fibered knots, the current authors improved this theorem as follows:

\begin{ther}{\rm\cite{Lagr}}\qua
{\rm{(a)}}\qua Let $X$ be any symplectic manifold with $b_2^+(X)>1$ which
contains an embedded symplectic torus with a vanishing cycle. Then for
each fibered knot $K$ in $S^3$, the result of knot surgery $X_K$ contains
infinitely many nullhomologous Lagrangian tori, pairwise inequivalent
under orientation-preserving diffeomorphisms.

\noindent{\rm{(b)}}\qua Let $X_i$, $i=1,2$, be symplectic 4--manifolds
containing embedded symplectic tori $F_i$ and assume that $F_1$ contains a
vanishing cycle. Let $X$ be the fiber sum, $X=X_1\#_{F_1=F_2}X_2$. Then
for each fibered knot $K$ in $S^3$, the manifold $X_K$ contains an
infinite family of homologically primitive and homologous Lagrangian
tori which are pairwise inequivalent.
\end{ther}

In sections \ref{LT}--\ref{L} we show how this theorem works in a
specific example constructed via double branched covers. The discussion
here differs somewhat from the more general arguments of \cite{Lagr},
however we feel that it is helpful to understand specific examples from
different points of view.

The authors gratefully acknowledge support from the National Science
Foundation. The first author was partially supported NSF Grants DMS0072212
and DMS0305818, and the second author by NSF Grant DMS0204041.

\section{Seiberg--Witten invariants}\label{sw}

The Seiberg--Witten invariant  of a smooth closed oriented 4--manifold $X$
with $b_2 ^+(X)>1$ is an integer-valued function which is defined on
the set of $\spin^{c}$ structures over $X$ (cf \cite{W}). In case
$H_1(X;\Z)$ has no 2--torsion, there is a natural identification of the
$\spin^{c}$ structures of $X$ with the characteristic elements of
$H_2(X;\Z)$ (ie those elements $k$ whose Poincar\'e duals $\hat{k}$
reduce mod~2 to $w_2(X)$).  In this case we view the Seiberg--Witten
invariant as
$$\ssw_X\co\lbrace k\in H_2(X;\Z)|\hat{k}\equiv
  w_2(TX)\pmod2)\rbrace \rightarrow \Z.$$
The sign of $\ssw_X$ depends on an orientation of $H^0(X;\R)\otimes\det
H_+^2(X;\R)\otimes \det H^1(X;\R)$; however, when $X$ has a symplectic
structure, there is a preferred sign for $\ssw_X$ (see \cite{T1}).

If $\ssw_X(\b)\neq 0$, then $\b$ is called a {\it{basic class}}
of $X$. It is a fundamental fact that the set of basic classes is
finite. Furthermore, if $\b$ is a basic class, then so is $-\b$ with
$\ssw_X(-\b)=(-1)^{(\text{e}+\text{sign})(X)/4}\,\ssw_X(\b)$ where
$\text{e}(X)$ is the Euler number and $\text{sign}(X)$ is the signature
of $X$.

It is convenient to view the Seiberg--Witten invariant as an element of
the integral group ring $\Z H_2(X)$. For $\a\in H_2(X)$ we let $t_\a$
denote the corresponding element in $\Z H_2(X)$. More specifically,
suppose that
$\{\pm \b_1,\dots,\pm \b_n\}$ is the set of nonzero basic classes for $X$.
Then the Seiberg--Witten invariant of $X$ is the Laurent polynomial
$$\sw_X = \ssw_X(0)+\sum_{j=1}^n \ssw_X(\b_j)\cdot
  (t_{\b_j}+(-1)^{(\text{e}+\text{sign})(X)/4}\, t_{\b_j}^{-1}) \in
  \Z H_2(X).$$
A key vanishing theorem for the Seiberg--Witten invariants is:

\begin{ther}{\rm\cite{W}}\qua
Let $X$ be a smooth closed 4--manifold which admits a decomposition
$X=A\cup B$ into 4--manifolds with $\bd A=\bd B = Y$. Suppose that
$b_2^+(A)>0$, $b_2^+(B)>0$, and that $Y$ admits a metric of positive
scalar curvature, then $\sw_X=0$.
\end{ther}

Another important and extremely useful fact about Seiberg--Witten
invariants is the adjunction inequality: If $X$ is a smooth closed
4--manifold with $b_2^+>1$ and $\Sig$ is an embedded surface of positive
genus $g$ in $X$ representing a nontrivial element of $H_2(X;\R)$ then
for any basic class $\b$ of $X$
\begin{equation}\label{adj}
2g-2\ge \Sig^2+\b\cdot\Sig
\end{equation}

We next recall the link surgery construction of \cite{KL4M}. This
construction starts with an oriented $n$--component link
$L=\{K_1,\dots,K_n\}$ in $S^3$ and $n$ pairs $(X_i,T_i)$ of smoothly
embedded self-intersection $0$ tori in simply connected 4--manifolds. (In
the original article \cite{KL4M}, an extra condition (`c-embedded')
was placed on these tori; however, recent work of Cliff Taubes \cite{T}
has shown this condition to be unnecessary.)

Let $\alpha_L\co\pi_1(S^3\setminus L)\to \Z $ denote the homomorphism
characterized by the property that it sends the meridian $m_i$ of
each component $K_i$ to $1$, and let $\ell_i$ denote the longitude of
$K_i$. The curves $\g_i=\ell_i - \alpha_L(\ell_i) m_i$ on $\bd N(K_i)$
form the boundary of a Seifert surface for the link, and in case $L$
is a fibered link, the $\gamma_i$ are given by the boundary components
of a fiber.

In $S^1\x (S^3\- N(L))$ let $T_{m_i}=S^1\x m_i$,
and define the link-surgery manifold $X(X_1,\dots X_n;L)$ by
$$X(X_1,\dots X_n;L)=
  (S^1\x (S^3\- N(L))\cup\bigcup\limits_{i=1}^n (X_i\- (T_i\x D^2))$$
where $S^1\x\bd N(K_i)$ is identified with $\bd N(T_i)$ so that for
each $i$
$$[T_{m_i}]=[T_i], \quad \text{and} \quad [\g_i] = [{\text{pt}}\x\bd D^2].$$
It is not clear whether or not this determines $X(X_1,\dots X_n;L)$
up to diffeomorphism, however any such manifold will have the same
Seiberg--Witten invariant:

\begin{ther}{\rm\cite{KL4M}}\qua
If each $\pi_1(X\setminus T_i) =1$, then $X(X_1,\dots X_n;L)$ is
simply-connected and its Seiberg--Witten invariant is
$$\sw_{X(X_1,\dots X_n;L)} =
  \Ds_L(t_1^2,\dots,t_n^2)\cdot\prod_{j=1}^n\sw_{X_j}\cdot(t_j-t_j^{-1})$$
where $t_j=t_{[T_j]}$ and $\Ds_L(t_1,\dots,t_n)$ is the symmetric
multivariable Alexander polynomial of the link $L$.
\end{ther}

In case each $(X_i,T_i)\cong (X,T)$, a fixed pair, we write \[ X(X_1,\dots
X_n;L)=X_L\] (We implicitly remember $T$, but it is removed from the
notation.)
As an example, consider the case where each $X_i=E(1)$, the rational elliptic
surface ($E(1)\cong \CP\#9\CPb$) and each $T_i=F$ is a smooth elliptic fiber.
Since $SW_{E(1)}=(t-t^{-1})^{-1}$, we have that
$$\sw_{E(1)_L}=\Ds_L(t_1^2,\dots,t_n^2).$$
In case the link $L$ is actually a knot $K$, we call the procedure
`knot surgery' and the resulting manifold $X_K$. The formula for the
Seiberg--Witten invariant looks slightly different in this case due to
the difference in the relationship of the Seiberg--Witten invariant of
a 3--manifold and its Alexander polynomial when $b_1>1$ and $b_1=1$.

\begin{ther}{\rm\cite{KL4M}}\qua
If $\pi_1(X\setminus T) =1$, then $X_K$ is
simply connected and its Seiberg--Witten invariant is
$$\sw_{X_K} =\Ds_K(t^2)\cdot\sw_X$$
where $t=t_{[T]}$.
\end{ther}

\section{Tori and simple covers}\label{sc}

Our first construction utilizes an extremely interesting theorem of
Jos\'e Montesinos and Hugh Morton which characterizes fibered links in the
3--sphere. To begin, let $X$ be a simply connected symplectic 4--manifold
containing an embedded symplectic torus $T$ of self-intersection $0$,
and identify a tubular neighborhood of $T$ with $S^1\x (S^1\x D^2)$. A
closed braid may be viewed as contained in $S^1\x D^2\C S^3 = (S^1\x
D^2)\cup (D^2\x S^1)$ and then its {\it axis} is $\{0\}\x S^1$. The
theorem of Montesinos and Morton is:

\begin{ther}[Montesinos and Morton \cite{MM}]
Every fibered link in $S^3$ with $k$ components can be obtained as the
preimage of the braid axis for a $d$--sheeted simple branched cover of
$S^3$ branched along a suitable closed braid, where $d=\max\{k,3\}$.
\end{ther}
\noindent (Recall that a {\it simple} branched cover of degree $d$
is one whose branch points have exactly $d-1$ points in their preimages.)

A second important ingredient in this construction is a theorem of
Kanenobu concerning the Hosokawa polynomial of fibered links. The
Alexander polynomial of a link $L$ of $k$ components is a polynomial
$\DD_L(t_1,\dots,t_k)$ in $k$ variables (corresponding to the meridians of
the components of the link). The polynomial $\DD_L(t,\dots,t)$ obtained
by setting all the variables equal is always divisible by $(t-1)^{k-2}$,
and the Hosokawa polynomial of $L$ is defined to be
$\DN_L(t)=\DD_L(t,\dots,t)/(t-1)^{k-2}$.

\begin{ther}[Kanenobu \cite{K}]
Let $f(t)$ be any symmetric polynomial of even degree with
integral coefficients satisfying $f(0)=\pm1$, then for any $k\ge2$
there is a fibered link $L$ of $k$ components in $S^3$ with
$\DN_L(t)=f(t)$.
\end{ther}

We now use these two theorems to build symplectic tori homologous to
multiples of $T$. We have described a tubular neighborhood $N$ of $T$
as $ N=S^1\x (S^1\x D^2)$. Fix a three-component fibered link $L$ in $S^3$
and let $B_L$ be the braid corresponding to $L$ by the Montesinos--Morton
Theorem. As above, we view $B_L$ as contained in $S^1\x D^2\C S^3 =
(S^1\x D^2)\cup (D^2\x S^1)$ with axis $A=\{0\}\x S^1$. (See Figure~1
for an example.) Then $T_L=S^1\x B_L\C N$ is a symplectic torus
\cite{surfaces}, and if $B_L$ has $m$ strands, then $T_L$ is homologous to $mT$.

\centerline{\unitlength 2cm
\begin{picture}(4.5,3)
\put(3,1.5){\oval(1.5,.8)[l]}
\put(3,1.5){\oval(1.5,.8)[t]}
\put(3.05,1.8){\line(0,-1){1}}
\put(3.25,1.8){\line(0,-1){1.2}}
\put(3.45,1.8){\line(0,-1){1.4}}
\curve(3.5,1.1,3.7,1.2,3.75,1.5)
\curve(3.1,1.1,3.2,1.1)
\curve(3.3,1.1,3.4,1.1)
\put(3.05,.8){\line(-1,0){.7}}
\put(3.25,.6){\line(-1,0){.9}}
\curve(2.35,.8,1.95,.4)
\curve(2.35,.6,2.27,.68)
\curve(2.15,.8,2.23,.72)
\put(3.45,.4){\line(-1,0){1.3}}
\curve(2.15,.4,2.07,.48)
\curve(1.95,.6,2.03,.52)
\put(1.95,.4){\line(-1,0){1}}
\put(1.95,.6){\line(-1,0){.8}}
\put(2.15,.8){\line(-1,0){.8}}
\put(1.35,.8){\line(0,1){1.5}}
\put(1.15,.6){\line(0,1){1.9}}
\put(.95,.4){\line(0,1){2.3}}
\put(1.35,2.3){\line(1,0){1.7}}
\put(1.15,2.5){\line(1,0){2.1}}
\put(.95,2.7){\line(1,0){2.5}}
\put(3.45,2.7){\line(0,-1){.7}}
\put(3.25,2.5){\line(0,-1){.5}}
\put(3.05,2.3){\line(0,-1){.3}}
\put(3.7,1.8){\small$A$}
\put(1.5,2.1){\small$B_L$}
\put(1.9,.1){\small Figure 1}
\end{picture}}

Let $\pi\co(S^3,L)\to (S^3,A)$ be the threefold branched cover with branch
set $B_L$ given by the Montesinos--Morton Theorem. Because $L=\pi^{-1}(A)$
is a three-component link, the covering restricted to $A$ is trivial. This
means that the restriction of $\pi$ over $\bd(S^1\x D^2)$ is a trivial
covering, and the induced branched cover over $N=S^1\x (S^1\x D^2)$
extends trivially over $X$. We thus get a threefold simple branched cover
$p=p_L\co\X\to X$ with branch set $T_L$. We have
$$\X = \widetilde{N}\cup\bigcup\limits_{i=1}^3(X\- N)_i$$
where $\widetilde{N}=p^{-1}(N)=S^1\x (S^3\- L)$, and $(X\- N)_i$ denotes a
copy of $X\- N$.

This means that $\X$ is obtained via link surgery on the link $L$ using
$(X,T)$. The Seiberg--Witten invariant of $\X$ (viewed as an element of
$\Z H_2(\X)$) may be calculated via the techniques of \cite{T, P, KL4M}:
$$\sw_{\X}=\Ds_L(t_1^2,t_2^2,t_3^2)\cdot
  \prod\limits_{i=1}^3\sw_{X_i}\cdot(t_i-t_i^{-1})$$
The induced map $p_*\co\Z H_2(\X)\to \Z H_2(X)$ satisfies
$p_*(\sw_{X_i})=\sw_X$. Also, since $t_i$ is the element of $\Z H_2(\X)$
corresponding to the homology class of $S^1\x \mu_i$
where $\mu_i$ are the meridians of the components of L, $p_*(t_i)$ is
the element of $\Z H_2(X)$ corresponding to $S^1\x\mu_A$, where $\mu_A$
is a meridian to $A$. Since $\mu_A $ is the core circle $S^1\x\{ 0\}\C
S^1\x D^2$, we have $[S^1\x\mu_A]=[T]$ in $H_2(X)$. Thus $p_*(t_i)=t$, and
\begin{equation}\label{sw1} 
p_*(\sw_{\X}) = \Ds_L(t^2,t^2,t^2)\cdot \sw_X^3\cdot(t-t^{-1})^3
\end{equation}
Now suppose that we are given another three-component link $L'$ which is a
threefold simple cover of $S^3$ with branch set $B_{L'}$ and symplectic
torus $T_{L'}=S^1\x B_{L'}$. The covering projections $p_L$, $p_{L'}$
are determined by homomorphisms $\vp_L$ (or $\vp_{L'})$ from $\pi_1(X\-
T_L)$ (or $\pi_1(X\- T_{L'})$) to the symmetric group $S_3$ such that
each meridian of $T_L$ (or $T_{L'}$) is sent to a transposition.

Any isotopy of $X$ taking $T_L$ to $T_{L'}$ and which carries the covering
data for $p_L$ to that of $p_{L'}$ gives rise to
\begin{equation}\label{cd1}
 \begin{CD} 
 {\X_L} @> {\widetilde{f}\  \cong} >> { \X_{L'}} \\ 
  @ V{p_L} VV  @ VV {p_{L'}} V \\ 
 {X} @>> {f\ \cong} > {X} \\ 
 \end{CD} 
\end{equation}
where $f(T_L)=T_{L'}$ and $f_*$ is the identity on homology.  

Since $\widetilde{f}_*(\sw_{\X_L})= \sw_{\X_{L'}}$ it follows from
\eqref{sw1} and \eqref{cd1} that
$$\Ds_L(t^2,t^2,t^2)=\Ds_{L'}(t^2,t^2,t^2)$$
In other words, 
$\DNS_L(t^2)=\DNS_{L'}(t^2)$. Using Kanenobu's theorem, one sees
that there are infinite families of fibered links $\{ L_i\}$ whose
$\DNS_{L_i}(t)$ are distinct and have arbitrary fixed even degree ($>
0$). The genus $g_L$ of the fibered link $L$ is half the degree of its
Hosokawa polynomial. (See, for instance, \cite{BZ}.)
Furthermore, the fiber of $L$ is the thrice-punctured surface which
is a simple threefold branched cover of $D^2$ (a normal fiber to $S^1\x
\{0\}$) with $m$ branch points. Thus the number of strands $m$ of $B_L$
is determined by $m=2g_L+4$.

This means that for any even $m\ge6$ we get an infinite family $\{ T_i\}$
of symplectic tori homologous to $mT$ with distinct threefold simple
branched covers. Note that each braided torus $T_L$ admits at most
finitely many simple threefold branched covers of $X$ with $T_L$ as branch
set, since there are finitely many distinct homomorphisms $\pi_1(X\-
T_L)\to S_3$. Thus we have:

\begin{thm}
Let $X$ be a simply connected symplectic 4--manifold containing an
embedded symplectic torus $T$ of self-intersection $0$. Then for each
even $m\ge6$ there are infinitely many pairwise nonsmoothly isotopic
embedded symplectic tori homologous to $mT$.\endproof
\end{thm} 

\section{Fiber sums}\label{FibS}

We begin this section with the same hypotheses as the last: We are given
a simply connected symplectic 4--manifold $X$ containing an embedded
symplectic torus $T$ of self-intersection $0$. The construction of
new symplectic tori is similar to that of the last section (and of
\cite{surfaces}). For each $m\ge 2$ consider closed braids $B$ with $m$
strands. Then the braided torus $T_B=S^1\x B$ is embedded in the tubular
neighborhood $S^1\x (S^1\x D^2)$ of $T=S^1\x S^1\x\{ 0\}$. Furthermore,
$T_B$ is symplectic and homologous to $mT$.

Suppose that $B$ and $B'$ are $m$--strand closed braids and that $T_{B'}$
is smoothly isotopic to $T_B$ in $X$. Then there is a diffeomorphism
$f\co X\to X$ satisfying: $f(T_B)=T_{B'}$, $f(\mu_B)=\mu_{B'}$, and
$f_*={\rm id}$ on $H_*(X)$.  (Here $\mu_B$ and $\mu_{B'}$ are meridians
to the braids; so they also may be viewed as meridians to the tori $T_B$
and $T_{B'}$.)

Our goal is to use relative Seiberg--Witten invariants $\sw_{(X,T_B)}$
to distinguish the tori $T_B$ up to isotopy. Let $E(1)$ denote the
rational elliptic surface. Because of the gluing theorems of \cite{T,P}
and the fact that the relative Seiberg--Witten invariant of $E(1)$ minus
a a smooth elliptic fiber is $\sw_{E(1)\- F}=1$ (see eg \cite{McT}),
the relative Seiberg--Witten invariant of $(X,T_B)$ may be expressed as
the absolute Seiberg--Witten invariant of the fiber sum of $X$ and $E(1)$
along $T_B$ and $F$:
$$\sw_{(X,T_B)}=\sw_{X\#_{T_B=F}E(1)}$$
Now write $N(T_B)$ for a tubular neighborhood of $T_B$ in $X$ and also
write$N(T)=S^1\x(S^1\x D^2)$, the original tubular neighborhood of
$T$. We have
\begin{equation}\label{union}
X\- N(T_B) = {\big{(}}X\- N(T){\big{)}} \cup \left(S^1\x \left((S^1\x
D^2)\- N(B)\right)\right )
\end{equation}
Let $L_B$ be the link in $S^3$ consisting of the closed braid $B$
together with its axis $A$.  If $\mu_A$ denotes a meridian to $A$, then
$T$ is homologous to $S^1\x\mu_A$. Let $t=t_T$ denote the corresponding
element in $\Z H_2(X)$.

We may now rewrite  \eqref{union} as
$$X\- N(T_B) =\left( X\- N(T)\right) \cup \left( S^1\x (S^3\- N(L_B))\right)$$
The manifold $ X\#_{T_B=F}E(1)$ is obtained from the same components as
link surgery using the link $L_B$ and the manifolds $(E(1),F)$ and
$(X,T)$; however the gluings are not necessarily those specified in
section \ref{sw}.  Since $E(1)$ has big diffeomorphism group with respect
to $F$ (see eg \cite{GS}), each diffeomorphism $\bd N(F)\to \bd N(F)$
extends to a self-diffeomorphism of $E(1)\- N(F)$; so the diffeomorphism
used to glue in $E(1)\- N(F)$ is inconsequential. However, it is useful
to demand that the fiber $F$ of $E(1)$ should be identified with $S^1\x
\lam_B$ where $\lam_B$ is the longitude of $B$ in $S^3$.

According to \cite{T,P}, $\sw_X\cdot (t-t^{-1})$ is the relative
Seiberg--Witten invariant of $(X,T)$, and by \cite{KL4M}, as described in
section \ref{sw}, the relative invariant of 
the manifold $S^1\x (S^3\- N(L_B))$ is $\Ds_{L_B}(t^2,\t^2)$.
Applying \cite{KL4M} and \cite{T}  we obtain:
$$\sw_{(X,T_B)}=\sw_{X\#_{T_B=F}E(1)} =
  \Ds_{L_B}(t^2,\t^2)\cdot\sw_X\cdot (t-t^{-1})$$
where $\t$ is the element of $\Z H_2(X)$ corresponding to    
$[S^1\x \mu_B]$. 
Since $[F]=[S^1\x \lam_B]=m [S^1\x\mu_A]=m[T]$. When applying this
formula, we need to remember that $t_T=t$ and $t_F=t^m$.
 
\begin{thm}
Let $X$ be a simply connected symplectic 4--manifold with $b_2^+ > 1$
containing an embedded symplectic torus $T$ of self-intersection $0$. For
a fixed integer $m\ge2$, let $B$ and $B'$ be closed $m$-strand braids
in $S^3$. Then $T_B$ and $T_{B'}$ are embedded symplectic tori in $X$
which are homologous to $mT$. If there is an isotopy of $X$ taking $T_B$
to $T_{B'}$, then $\Ds_{L_{B'}}(t^2,\t'^2)=\Ds_{L_B}(t^2,\t^{\pm2})$.
\end{thm}   

\begin{proof}
We first describe $H_2(X\#_{T_B=F}E(1))$. Let $R_B$ denote the group
of rim tori of the torus $T_B$; ; ie $R_B=\ker(H_2(X\- T_B)\to
H_2(X))\cong \Z\oplus\Z$. A basis for $R_B$ is given by $\t=[S^1\x
\mu_B]$ and $v=[\lam_B\x \mu_B]$ where $\lam_B$ is the longitude of the
knot $B$ in $S^3$. The classes $\t$ and $v$ are primitive (because of
the definition of $R_B$), thus there is a group $D_B\cong\Z\oplus\Z$
generated by the dual classes to $\t$ and $v$ in $H_2(X\#_{T_B=F}E(1))$.

Let $A=T_B^{\perp}=T^{\perp}\C H_2(X)$. Note that the adjunction
inequality \eqref{adj} implies that no basic class of $X$ has nontrivial
intersection with $[T]$.
Thus $\sw_X\in\Z{A}$.  We have $H_2(X\- T_B)= A\oplus R_B$. Finally,
suppose that $[T]$ is $n$ times a primitive homology class, and let $S$
denote the class in $H_2(X\#_{T_B=F}E(1))$ which has a representative
built from $mn$ punctured sections in $E(1)\- F$ and a surface in $X\-
T_B$ which has boundary $mn$ copies of the meridian $\mu_B$ to $T_B$.

A Mayer--Vietoris argument shows that the homology of $X\#_{T_B=F}E(1)$
splits as
$$H_2(X\#_{T_B=F}E(1))=A\oplus (R_B \oplus D_B) \oplus \Z(S) \oplus E_8$$
where the $E_8$ comes from $H_2(E(1)\- F)$. There is a similar splitting
of the homology of $H_2(X\#_{T_{B'}=F}E(1))$.

If there is an isotopy of $T_B$ to $T_{B'}$, there is a diffeomorphism
$$\f\co X\#_{T_B=F}E(1)\to X\#_{T_{B'}=F}E(1)$$
satisfying $\f_*|_A={\rm id}$ \ and $\f_*(R_B)=R_{B'}$ (since
$f(\mu_B)=\mu_{B'}$). Thus the induced homomorphism of group rings
satisfies $\f_*(\sw_X)=\sw_X$ and $\f_*(t_F)=t_F$; ie $\f_*(t)^m=t^m$,
and so $\f_*(t)=t$ because $H_2(X\#_{T_B=F}E(1))$ is torsion-free.
It follows that the fact that $\f_*(\sw_{(X,T_B)})=\sw_{(X,T_{B'})}$
implies that
\begin{equation}\label{Delta}
\Ds_{L_B}(t^2,\f_*(\t)^2)=\Ds_{L_{B'}}(t^2,\t'^2)
\end{equation}
Write $\f_*(\t)=a\t'+bv'$ (where $\t'=[S^1\x \mu_{B'}]$ and
$v'=[\lam_{B'}\x \mu_{B'}]$). Each term $nt^{2r'}\t'^{2s'}$
of $\Ds_{L_{B'}}(t^2,\t'^2)$ corresponds to basic classes of
$X\#_{T_{B'}=F}E(1)$ of the form
$\a+(2r'\pm 1)[T]+2s'\t'$ where $\a\in H_2(X)$ is a basic class, and so 
$\a\cdot \t' =0$, and $\a\cdot [T]=0$. Furthermore each class in $R_{B'}$
is orthogonal to $T$.

Terms of the form $nt^{2r}\f_*(\t)^{2s}$ of $\Ds_{L_{B'}}(t^2,\f_*(\t)^2)$
correspond to basic classes of $X\#_{T_{B'}=F}E(1)$ of the form
$\b+(2r\pm 1)[T]+2s(a\t'+bv')$, and each basic class can be written like
this. Since $\t'$ and $v'$ are independent, it is clear that $b=0$. This
means that $\f_*(\t)=a\t'$, and $a=\pm1$ since $\t$ is primitive. Thus
$\Ds_{L_{B'}}(t^2,\t'^2)=\Ds_{L_B}(t^2,\t^{\pm2})$.
\end{proof}

We have as a corollary:

\begin{thm}
Let $X$ be a simply-connected symplectic 4--manifold satisfying
$b_2^+(X)>1$ and containing an embedded symplectic torus $T$
of self-intersection $0$. For each $m\ge2$ there are infinitely many
pairwise nonisotopic embedded symplectic tori in $X$ which are homologous
to $mT$. \end{thm}

\begin{proof}
This follows from the above theorem provided for each $m\ge 2$ there are
infinitely many closed m-strand braids $B$ whose 2-component links $L_B=
A \cup B$ have distinct 2-variable Alexander polynomials. Such examples
are given, for example, in the work of Etgu and Park \cite{EP}.
\end{proof}

\section{Lagrangian tori}\label{LT}

In this section we use branched covers as a means for constructing
examples of Lagrangian tori in symplectic 4--manifolds whose homology
classes are equal but which are not equivalent under symplectic
diffeomorphisms. There are already two papers \cite{V, Lagr} dealing
with this phenomenon, and the invariants of \cite{Lagr} can be used
to distinguish the examples given in this section. However we believe
that the constructions below are interesting in their own right and are
certainly different from those cited.

To begin, let $K$ be the trefoil knot, and $M_K$ the 3--manifold
obtained from 0--framed surgery on $S^3$ along $K$. Since $K$ is
a genus--1 fibered knot, $M_K$ fibers over the circle with fiber a
torus, $M_K=T^2\x_\vp S^1$. Let $E(1)_K$ be the result of knot surgery
on $E(1)$, $E(1)_K=E(1)\#_{F=S^1\x m_0}S^1\x M_K$, where $m_0$ is a
meridian to $K$. This manifold has a symplectic structure induced from
that on $E(1)$ and the structure on $S^1\x M_K$ in which the fiber and
section are symplectic submanifolds. (See \cite{KL4M}.)

\centerline{\unitlength 1cm
\begin{picture}(5.5,5.5)
\put (3,3){\oval(3,4)[r]}
\put (2,5){\line(3,-4){.4}}
\put (3,3.66){\line(-3,-4){1}}
\put (3,3.66){\line(-3,4){.4}}
\put (2,3.66){\line(3,4){1}}
\put (2,3.66){\line(3,-4){.4}}
\put (3,2.33){\line(-3,4){.4}}
\put (2,2.33){\line(3,-4){.4}}
\put (3,2.33){\line(-3,-4){1}}
\put (3,1){\line(-3,4){.4}}
\curve(2,5,1,5,0.6,4.75,0.5,4)
\curve(.4,4.3, .35,4.29, .3,4.27, .25,4.22, .225,4.15, .224,4.05, .275,3.95, .35,3.92, .4,3.915, .52,3.915, .57,3.92, .645, 3.95, .696,4.05, .695,4.15, .67,4.22,.62,4.27, .57,4.3 )
\curve(.4,3.27, .35,3.26, .3,3.24, .25,3.19, .225,3.12, .224,3.02, .275,2.92, .35,2.89, .4,2.885, .52,2.885, .57,2.89, .645, 2.92, .696,3.02, .695,3.12, .67,3.19,.62,3.24, .57,3.27 )
\curve(.4,2.25, .35,2.24, .3,2.22, .25,2.17, .225,2.1, .224,2, .275,1.9, .35,1.87, .4,1.865, .52,1.865, .57,1.87, .645, 1.9, .696,2, .695,2.1, .67,2.17,.62,2.22, .57,2.25 )
\curve(2,1,1,1,0.6,1.25,.5,1.75)
\put(-.25,3.7){\small{$m_1$}}
\put(-.25,1.7){\small{$m_2$}}
\put(-.25,2.72){\small{$m_0$}}
\put(0,1.0){$K$}
\curve(.5,3.8,.5,3.02)
\curve(.5,2.8,.5,2)
\put(1.85,.4){\small Figure 2}
\end{picture}}

Let $m_1$ and $m_2$ be meridians of $K$ as in Figure~2, and let $X$
denote the double branched cover of $E(1)_K$ with branch set $S^1\x
(m_1\cup m_2)$. Since $S^1\x m_i$ is a section to the fibration $S^1\x
M_K\to S^1\x S^1$, the branch set of this cover is symplectic, hence $X$
inherits a symplectic structure. We have
$$X= E(1)'\#_{F'=S^1\x m'_0} \, S^1\x \widetilde{M}_K\#_{S^1\x
  m''_0=F''}\, E(1)''$$
where $E(1)'$ and $E(1)''$ are copies of $E(1)$ and $ \widetilde{M}_K$
is the double cover of $M_K$ branched over $m_1\cup m_2$. It follows that
$\widetilde{M}_K$ also fibers over the circle, and its fiber is the double
branched cover of the fiber of $M_K$, branched over two points. Thus
the fiber of $\widetilde{M}_K\to S^1$ has genus 2. We can say more:

\begin{lem}
Let $K$ be any knot in $S^3$ and $M_K$ the result of 0--surgery along
$K$. The double cover of $M_K$ branched over two meridians to $K$ is
$M_{K\#  K}$, the result of 0--surgery on $S^3$ along the connected
sum $K\# K$.
\end{lem}

\begin{proof}
This proof is an exercise in Kirby calculus. The double branched cover of
$S^3$ branched over the two-component unlink is $S^2\x S^1$. This means
that the double branched cover of $M_K$ branched along two meridians
to $K$ is the result of surgery on the lift of $K$ in $S^2\x S^1$. (See
Figure~3.) Note that $K$ lifts to two components. Referring to Figure~3,
slide one copy of  $K$ over the other copy of $K$ to obtain Figure~4. In
this figure, 0--surgery on $K$ together with 0--surgery on a meridian
form a cancelling pair. We are left with 0--surgery on $K\#  K$.
\end{proof}

\centerline{\unitlength .8cm
\begin{picture}(7,5.5)
\curve(1.9,1.45,1.75,1.6,1.5,1.9,1.25,2.4,1.22,2.7,1.25,3,1.5,3.5,1.75,3.8,2,4,2.25,4.1,2.75,4.2,3,4.22,3.5,4.25,4,4.25,4.5,4.22,4.75,4.2,5.25,4.1,5.5,4,5.75,3.8,6,3.5,6.25,3,6.28,2.7,6.25,2.4,6,1.9,5.75,1.6,5.6,1.45)
\curve(2,4.2,1.9,4.5,1.8,4.7,1.6,5,1.4,5.2,1.2,5.3,1,5.35,.8,5.35,.6,5.3,.4,5.2,.2,5.05,0,4.7,-.1,4.5,-.2,4.2,-.25,3.8,-.3,3.3,-.31,3.15)
\curve(-.31,2.25,-.3,2.1,-.25,1.6,-.2,1.2,
-.1,.9,0,.7,.2,.35,.4,.2,.6,.1,.8,.05,1,.05,1.2,.1,1.4,.2,1.6,.4,1.8,.7,1.9,.9,2,1.2,2.05,1.6,2.1,
2.1,2.13,2.7,2.1,3.3,2.05,3.8)
\curve(2.25,1.3,2.75,1.2,3,1.18,3.5,1.15,4,1.15,4.5,1.18,4.75,1.2,5.25,1.3)
\curve(5.5,4.2,5.6,4.5,5.7,4.7,5.9,5,6.1,5.2,6.3,5.3,6.5,5.35,6.7,5.35,6.9,5.3,7.1,5.2,7.3,5.05,7.5,4.7,7.6,4.5,7.7,4.2,7.75,3.8,7.8,3.3,7.81,3.15)
\curve(7.81,2.25,7.8,2.1,7.75,1.6,7.7,1.2,
7.6,.9,7.5,.7,7.3,.35,7.1,.2,6.9,.1,6.7,.05,6.5,.05,6.3,.1,6.1,.2,5.9,.4,5.7,.7,5.6,.9,5.5,1.2,5.45,1.6,5.4,2.1,5.37,2.7,5.4,3.3,5.45,3.8)
\put(-.73,2.3){\framebox(.8,.8){$K$}}
\put(7.43,2.3){\framebox(.8,.8){$K$}}
\put(0,5.2){$0$}
\put(7.3,5.2){$0$}
\put(3.75,4.5){$0$}
\put(3,2.6){\vector(1,0){1.5}}
\put (2.85,-.25){\small Figure 3}
\end{picture}}
 
\bigskip

\centerline{\unitlength .8cm
\begin{picture}(16,5.75)
\curve(4.13,2.9,4.1,3.3,4.05,3.8,4,4.2,3.9,4.5,3.8,4.7,3.6,5,3.4,5.2,3.2,5.3,3,5.35,2.8,5.35,2.6,5.3,2.4,5.2,2.2,5.05,2,4.7,1.9,4.5,1.8,4.2,1.75,3.8,1.7,3.3,1.69,3.15)
\curve(1.69,2.25,1.7,2.1,1.75,1.6,1.8,1.2,
1.9,.9,2,.7,2.2,.35,2.4,.2,2.6,.1,2.8,.05,3,.05,3.2,.1,3.4,.2,3.6,.4,3.8,.7,3.9,.9,4,1.2,4.05,1.6,4.1,2.1,4.13,2.5)
\curve(7.81,2.25,7.8,2.1,7.75,1.6,7.7,1.2,
7.6,.9,7.5,.7,7.3,.35,7.1,.2,6.9,.1,6.7,.05,6.5,.05,6.3,.1,6.1,.2,5.9,.4,5.7,.7,5.6,.9,5.5,1.2,5.45,1.6,5.4,2.1,5.37,2.5)
\curve(5.37,2.9,5.4,3.3,5.45,3.8,5.5,4.2,5.6,4.5,5.7,4.7,5.9,5,6.1,5.2,6.3,5.3,6.5,5.35,6.7,5.35,6.9,5.3,7.1,5.2,7.3,5.05,7.5,4.7,7.6,4.5,7.7,4.2,7.75,3.8,7.8,3.3,7.81,3.15)
\curve(4.13,2.9,5.37,2.9)
\curve(4.13,2.5,5.37,2.5)
\put(1.27,2.3){\framebox(.8,.8){$K$}}
\put(7.43,2.3){\framebox(.8,.8){$K$}}
\put(7.3,5.2){$0$}
\curve(13.25,1.3,12.75,1.2,12.5,1.18,12,1.15,11.5,1.15,11,1.18,10.75,1.2,10.25,1.3,9.9,1.45,9.75,1.6,9.5,1.9,9.25,2.4,9.22,2.7,9.25,3,9.5,3.5,9.75,3.8,10,4,10.25,4.1,10.75,4.2,11,4.22,11.5,4.25,12,4.25,12.5,4.22,12.75,4.2,13.25,4.1,13.5,4,13.75,3.8,14,3.5,14.25,3,14.28,2.7,14.25,2.4,14,1.9,13.75,1.6,13.6,1.45)
\curve(15.81,2.25,15.8,2.1,15.75,1.6,15.7,1.2,
15.6,.9,15.5,.7,15.3,.35,15.1,.2,14.9,.1,14.7,.05,14.5,.05,14.3,.1,14.1,.2,13.9,.4,13.7,.7,13.6,.9,13.5,1.2,13.45,1.6,13.4,2.1,13.37,2.7,13.4,3.3,13.45,3.8)
\curve(13.5,4.2,13.6,4.5,13.7,4.7,13.9,5,14.1,5.2,14.3,5.3,14.5,5.35,
14.7,5.35,14.9,5.3,15.1,5.2,15.3,5.05,15.5,4.7,15.6,4.5,15.7,4.2,
15.75,3.8,15.8,3.3,15.81,3.15)
\put(15.43,2.3){\framebox(.8,.8){$K$}}
\put(11.75,4.5){$0$}
\put(15.3,5.2){$0$}
\put (7.5,-.25){\small Figure 4}
\end{picture}}

\medskip
Since $X= E(1)'\#_{F'=S^1\x m'_0} \, S^1\x M_{K\# K}\#_{S^1\x m''_0=F''}\,
E(1)''$, and the complement of a fiber in $E(1)$ is simply-connected,
we have
$$\pi_1(X)= \pi_1(S^1\x m'_0\x\mu')\backslash\pi_1(S^1\x (M_{K\# K}\-
  (m_0'\cup m_0''))/\pi_1(S^1\x m''_0\x\mu'')$$
where $\mu'$ and $\mu''$ are the meridians to $m_0'$ and $m_0''$.
The group in the middle, $\pi_1(S^1\x M_{K\# K}\- (m_0'\cup m_0''))$,
is normally generated by the classes of $S^1\x\pt$, any meridian to $K\#
K$, and by $\mu'$ and $\mu''$. These loops all lie on $S^1\x m'_0\x\mu'$
or $S^1\x m''_0\x\mu''$; so we see that $X$ is simply-connected.

Consider the paths $P$ and $P'$ shown in Figure~5, each running from a
point $y_1\in m_1$ to $y_2\in m_2$. These paths lie in a fiber of the
fibration of $S^1\x M_K$ over $S^1\x S^1$; thus the construction of the
symplectic structure (essentially `area form on base' plus `area form on
fiber') implies that the surfaces $S^1\x P$ and $S^1\x P'$ are Lagrangian.

\centerline{\unitlength 1.1cm
\begin{picture}(4.5,5)
\put (3,3){\oval(3,4)[r]}
\put (3,3){\oval(3,4)[r]}
\put (2,5){\line(3,-4){.4}}
\put (3,3.66){\line(-3,-4){1}}
\put (3,3.66){\line(-3,4){.4}}
\put (2,3.66){\line(3,4){1}}
\put (2,3.66){\line(3,-4){.4}}
\put (3,2.33){\line(-3,4){.4}}
\put (2,2.33){\line(3,-4){.4}}
\put (3,2.33){\line(-3,-4){1}}
\put (3,1){\line(-3,4){.4}}
\curve(2,5,1,5,0.6,4.75,0.5,4)
\curve(.4,4.3, .35,4.29, .3,4.27, .25,4.22, .225,4.15, .224,4.05, .275,3.95, .35,3.92, .4,3.915, .52,3.915, .57,3.92, .645, 3.95, .696,4.05, .695,4.15, .67,4.22,.62,4.27, .57,4.3 )
\curve(.4,3.27, .35,3.26, .3,3.24, .25,3.19, .225,3.12, .224,3.02, .275,2.92, .35,2.89, .4,2.885, .52,2.885, .57,2.89, .645, 2.92, .696,3.02, .695,3.12, .67,3.19,.62,3.24, .57,3.27 )
\curve(.4,2.25, .35,2.24, .3,2.22, .25,2.17, .225,2.1, .224,2, .275,1.9, .35,1.87, .4,1.865, .52,1.865, .57,1.87, .645, 1.9, .696,2, .695,2.1, .67,2.17,.62,2.22, .57,2.25 )
\curve(2,1,1,1,0.6,1.25,.5,1.75)
\put(-.25,3.7){\small{$m_1$}}
\put(-.25,1.7){\small{$m_2$}}
\put(-.25,2.72){\small{$m_0$}}
\curve(.5,3.8,.5,3.02)
\curve(.5,2.8,.5,2)
\put (1.85,.5){\small Figure 5}
\curve(.645,3.95,1.1,3.1,.65,2.18)
\curve(.645,3.95,1.5,4.4,2.4, 4.37)
\curve(2.55, 4.32,3.2, 3.67,2.55, 3.02)
\curve(.65,2.18,2.35, 2.94)
\put(1.17,3.1){\small$P$}
\put(3.27,3.67){\small$P'$}
\put(1.1,3.1){\vector(0,-1){.01}}
\put(3.2,3.67){\vector(0,-1){.01}}
\end{picture}}

Since the endpoints of $P$ and $P'$ lie in the branch set of the cover,
their lifts $\g$ and $\g'$ in $\widetilde{M}_K=M_{K\# K}$ are circles in
the fibers (which are genus--2 surfaces). We thus obtain Lagrangian tori
$T=T_\g=S^1\x\g$ and $T'=T'_\g=S^1\x\g'$ in $S^1\x M_{K\# K}$. These tori
are disjoint from the lifts of $m_0$, where the gluing in the construction
of $X$ takes place, so $T$ and $T'$ are Lagrangian tori in $X$.

\centerline{\unitlength 1.2cm
\begin{picture}(8,6.15)
\put (2,5){\line(3,-4){.4}}
\put (3,3.66){\line(-3,4){.4}}
\put (2,3.66){\line(3,4){1}}
\put (2,3.66){\line(3,-4){.4}}
\put (2,2.33){\line(3,4){1}}
\put(3,2.33){\line(-3,4){.4}}
\put (2,2.33){\line(3,-4){.4}}
\put(3,1){\line(-3,4){.4}}
\put (2,1){\line(3,4){1}}
\curve(2,5,1,5,.6,4.75,.5,4)
\curve(.4,3.27, .35,3.26, .3,3.24, .25,3.19, .225,3.12, .224,3.02, .275,2.92, .35,2.89, .4,2.885, .52,2.885, .57,2.89, .645, 2.92, .696,3.02, .695,3.12, .67,3.19,.62,3.24, .57,3.27 )
\curve(2,1,1,1,0.6,1.25,.5,1.75)
\put(-.25,2.72){\small{$m_0''$}}
\curve(.5,4,.5,3.02)
\curve(.5,2.8,.5,1.75)
\curve(2.4,4.3,1.75,3.65,2.4,3)
\curve(2.6,4.3,3.25,3.65,2.6,3)
\put(1.75,3.65){\vector(0,-1){.01}}
\put(1.6,4){\small{$\delta''$}}
\put(3,5){\line(1,0){1.1}}
\put(3,1){\line(1,0){1.1}}
\curve(3.8,5.1,4,5.5,4.2,5,4,4.5,3.8,4.9)
\curve(3.8,1.1,4,1.5,4.2,1,4,.5,3.8,.9)
\put(3.3,5.3){\small$\widetilde{m}_1$}
\put(3.3,1.3){\small$\widetilde{m}_2$}
\put(4.3,5){\line(1,0){.7}}
\put(4.3,1){\line(1,0){.7}}
\put (5,5){\line(3,-4){.4}}
\put (6,3.66){\line(-3,-4){1}}
\put (6,3.66){\line(-3,4){.4}}
\put (6,5){\line(-3,-4){1}}
\put (5,3.66){\line(3,-4){.4}}
\put (5,2.33){\line(3,4){.4}}
\put (5,2.33){\line(3,-4){.4}}
\put (6,2.33){\line(-3,4){.4}}
\put (6,2.33){\line(-3,-4){1}}
\put (5,1){\line(3,4){.4}}
\put (6,1){\line(-3,4){.4}}
\curve(6,5,7,5,7.4,4.75,7.5,4)
\curve(7.5,4,7.5,3.02)
\curve(6,1,7,1,7.4,1.25,7.5,1.75)
\curve(7.5,2.8,7.5,1.75)
\curve(5.4,4.3,4.75,3.65,5.4,3)
\curve(5.6,4.3,6.25,3.65,5.6,3)
\put(6.245,3.65){\vector(0,-1){.01}}
\put(6.2,4){\small{$\delta'$}}
\curve(7.6,3.27, 7.65,3.26,7.7,3.24,7.75,3.19, 7.775,3.12, 7.776,3.02, 7.725,2.92, 7.65,2.89, 7.6,2.885,7.48,2.885,7.43,2.89,7.355, 2.92, 7.304,3.02, 7.305,3.12, 7.33,3.19,7.38,3.24, 7.43,3.27 )
\put(7.75,2.72){\small{$m_0'$}}
\put(3.3,0){\small Figure 6}
\end{picture}}

\medskip
The meridian $m_0$ in $M_K$ lifts to a pair of meridians, $m_0'$, $m_0''$
in $M_{K\# K}$ as in Figure~6. $H_1( M_{K\# K}\- (m_0'\cup m_0''))\cong
\Z\oplus\Z$ is generated by $[m_0']=[m_0'']$ and by the classes of the
meridians $[\mu']=[\mu'']$ to $m_0'$ and $m_0''$. (Note that $M_{K\# K}\-
(m_0'\cup m_0'')$ is fibered over the circle and its fibers are genus 2
surfaces with two boundary components. The meridians $\mu'$ and $\mu''$
form the boundary of one fiber.)

Referring to Figure~6, in $H_1( M_{K\# K}\- (m_0'\cup m_0''))$ we have
$[\g']-[\g]=[\d']+[\d'']$. Because $\d'$ and $\d''$ link neither $m_0'$
(or $m_0''$) nor $\mu'$ (or $\mu''$), the loops $\d'$ and $\d''$ are
nullhomologous in $M_{K\# K}\- (m_0'\cup m_0'')$. This means that the
corresponding Lagrangian tori, $\Sig'=S^1\x \d'$ and $\Sig''=S^1\x\d''$
are nullhomologous in $X$. Since $T'-T$ is homologous to $\Sig'+\Sig''$,
we see that $T$ and $T'$ are homologous in $X$.

The loop $\g$ is a separating curve in the fiber of $M_{K\# K}\to
S^1$. See Figure~7. We see that $\g$ is homologous to $\mu'$ in the
fiber of $M_{K\# K}\- (m_0'\cup m_0'')\to S^1$.  Thus in $S^1\x M_{K\#
K}\- (m_0'\cup m_0'')$, the Lagrangian torus $T=S^1\x \g$ is homologous
to $R_{\mu'}=S^1\x\mu'$. But $R_{\mu'}$ is a rim torus to $S^1\x m_0'$,
one of the tori along which the fiber sum
$$X= E(1)'\#_{F'=S^1\x m'_0} \, S^1\x M_{K\# K}\#_{S^1\x m'_0=F''}\, E(1)''$$
is made. In such a fiber sum, the rim tori give essential homology
classes -- thus we see that the Lagrangian tori $T$ and $T'$ are
essential in $X$.

\centerline{\unitlength 1.5cm
\begin{picture}(9,4)
\qbezier(2.75,3)(4.5,4)(6.25,3)
\qbezier(2.75,1)(4.5,0)(6.25,1)
\qbezier(2.75,3)(1,2)(2.75,1)
\qbezier(6.25,3)(8,2)(6.25,1)
\qbezier(2.85,2.1)(3.25,1.35)(3.65,2.1)
\qbezier(2.9,2)(3.25,2.6)(3.6,2)
\qbezier(5.35,2.1)(5.75,1.35)(6.15,2.1)
\qbezier(5.4,2)(5.75,2.6)(6.1,2)
\curve(4.5,3.5,4.9,2,4.5,.5)
\curve[20](4.5,3.5,4.15,2,4.5,.5)
\curve(1.9, 2.2,2, 2,1.9, 1.8)
\curve[10](1.9, 2.2,1.92, 2,1.9, 1.8)
\put(2.1,2.1){\small$\mu''$}
\curve(7.1, 2.2,7, 2,7.1, 1.8)
\curve[10](7.1, 2.2,7.08, 2,7.1, 1.8)
\put(6.75,2.1){\small$\mu'$}
\curve(4.5,3.5,6.4, 2,4.5,.5)
\curve[30](4.5,3.5,2.6, 2,4.5,.5)
\put(5.5, 2.67){\small$\g'$}
\put(4.9,2.4){\small$\g$}
\put(4.15,0) {\small Figure 7}
\end{picture}}

\smallskip 
We claim that there is no diffeomorphism of $X$ which takes $T$ to
$T'$. To see this we shall use an invariant obtained from Seiberg--Witten
theory. To do this we need the notion of `surgery on $T$'.  As usual,
this means the result of removing a tubular neighborhood $N(T)\cong
T^2\x D^2$ and regluing it.
$$X(T,\psi) = (X\- N(T)) \cup_\psi (T^2\x D^2)$$
The key quantity in this operation is the class $\o\in H_1(\bd N(T))$
which is killed by the composition of $\psi\co\bd N(T) \to T^2\x \bd D^2$
with the inclusion $T^2\x \bd D^2\to T^2\x D^2$. This class determines
$X(T,\psi)$ up to diffeomorphism; so we write $X_T(\o)$ instead of
$X(T,\psi)$.

Note that if there  is a diffeomorphism $f$ of $X$ taking $T$
to $T'$, then each manifold $X_T(\o)$ corresponds to a unique
$X_{T'}(f_*(\o))$. Thus the collection of all manifolds
$$\{ X_T(\o) |\,\o\in H_2(\bd N(T))\}$$
is a diffeomorphism invariant of $(X,T)$. Our invariant $I(X,T)$, defined
and computed below, will be the set of Seiberg--Witten invariants of
these manifolds.

\section{Product formulas for the Seiberg--Witten invariant}

Before formally defining $I(X,T)$ we need to discuss techniques for
calculating the Seiberg--Witten invariants of the manifolds $X_T(\o)$.
Fix simple loops $\a$, $\b$, $\d$  on $\bd N(T)$ whose homology classes
generate $H_1(\bd N(T))$. If $\o = p\a +q\b + r\d$ write $X_T(p,q,r)$
instead of $X_T(\o)$.  An important formula for calculating the
Seiberg--Witten invariants of surgeries on tori is due to Morgan,
Mrowka, and Szabo \cite{MMS} (see also \cite{MT}, \cite{T}). Given a
class $k\in H_2(X)$:
\begin{multline}\label{surgery formula}
\sum_i\ssw_{X_T(p,q,r)}(k_{(p,q,r)}+2i[T])=
  p\sum_i\ssw_{X_T(1,0,0)}(k_{(1,0,0)}+2i[T]) \\
+q\sum_i\ssw_{X_T(0,1,0)}(k_{(0,1,0)}+2i[T])
  +r\sum_i\ssw_{X_T(0,0,1)}(k_{(0,0,1)}+2i[T])
\end{multline}
In this formula, $T$ denotes the torus which is the core $T^2\x {0}\C
T^2\x D^2$ in each specific manifold $X(a,b,c)$ in the formula, and
$k_{(a,b,c)}\in H_2(X_T(a,b,c))$ is any class which agrees with the
restriction of $k$ in $H_2(X\- T\x D^2,\bd)$ in the diagram:
$$\begin{array}{ccc}
H_2(X_T(a,b,c)) &\longrightarrow & H_2(X_T(a,b,c), T\x D^2)\\
&&\Big\downarrow \cong\\
&&H_2(X\- T\x D^2,\bd)\\
&&\Big\uparrow \cong\\
H_2(X)&\longrightarrow & H_2(X,T\x D^2)
\end{array}$$
Let $\pi(a,b,c)\co H_2(X_T(a,b,c))\to H_2(X\- T\x D^2,\bd)$ be the
composition of maps in the above diagram, and $\pi(a,b,c)_*$
the induced map of integral group rings. Since we are
interested in invariants of the pair $(X,T)$, we shall work with
$$\swb_{(X_T(a,b,c),T)}=\pi(a,b,c)_*(\sw_{X_T(a,b,c)})\in \Z H_2(X\-
  T\x D^2,\bd).$$
The indeterminacy due to the sum in \eqref{surgery formula} is caused
by multiples of $[T]$; so passing to $\swb$ removes this indeterminacy,
and the Morgan--Mrowka--Szabo formula becomes
\begin{equation}\label{surgery formula 2}
 \swb_{(X_T(p,q,r),T)} = 
p\swb_{(X_T(1,0,0),T)} +q\swb_{(X_T(0,1,0),T)}
+r\swb_{(X_T(0,0,1),T)}.
\end{equation}

\begin{prop}
The collection of Seiberg--Witten invariants 
$$I(X,T)=\{ \swb_{X_T(a,b,c)}|\, a,b,c\in\Z\}$$
is an orientation-preserving diffeomorphism invariant of the pair
$(X,T)$.\endproof
\end{prop}

\section{Calculation of $I(X,T)$: $X(0,1,0)$}

We first specify a basis for $H_1(\bd N(T))$ as follows: Recall that
$T=S^1\x\g$ where $\g$ lies in a fiber of the fibration $M_{K\# K}\-
(m_0'\cup m_0'')\to S^1$ (Figure~7). Then $N(T)$ may be identified with
$S^1\x\g\x D^2$, and we take the basis $\a=[S^1\x\pt\x\pt]$, $\b=
[\pt\x\g_L]$, where $\g_L$ is a pushoff of $\g$ in the fiber of the
fibration $M_{K\# K}\- (m_0'\cup m_0'')\to S^1$ (this is called
the `Lagrangian framing' in \cite{Lagr}), and $\d=[m_T]$, where
$m_T=\pt\x\pt\x\bd D^2$, the meridian to $T$. It is then clear from
\eqref{surgery formula 2} that in order to calculate $I(X,T)$, one needs
to calculate $\swb_Y$ for $Y=X_T(1,0,0)$, $X_T(0,1,0)$, and $X_T(0,0,1)$;
however, from our choice of basis, we have $X_T(0,0,1)\cong X$.
This leaves us with two invariants to calculate below. For a different
approach to these calculations see \cite{Lagr}.

The calculation of $\sw_{X_T(0,1,0)}$ depends on some basic facts about
double covers of 3--manifolds branched over closed braids. Suppose that
$B$ is a braid in a solid torus with $2m$ strands, ie $B\C S^1\x D^2\C
M^3$ is a link such that each disk $\pt\x D^2$ intersects $B$ in exactly
$2m$ points. There is then a double cover $Y_B\to M^3$ branched over $B$
for which each meridian to $B$ is covered nontrivially, and this cover
is trivial outside the solid torus.

The pertinent question is: `What is the effect on $Y_B$ of putting
half-twists into the braid $B$?' In other words, suppose that $\z$ is an
arc in $(\pt\x D^2)\C S^1\x D^2$ whose endpoints lie on $B$, but so that
no other point of $\z$ is on $B$. Then we can put half-twists in $B$ by
twisting in a small neighborhood of $\z$. Figure~8 shows a local picture.

\centerline{\unitlength 1cm
\begin{picture}(8,2.25)
\put(0,.8){\line(1,0){3}}
\put(0,1.8){\line(1,0){3}}
\put(1.5,.8){\line(0,1){1}}
\put(1.2,1.2){\small{$\z$}}
\put(5,.8){\line(1,0){1}}
\put(5,1.8){\line(1,0){1}}
\put(7,.8){\line(1,0){1}}
\put(7,1.8){\line(1,0){1}}
\put(6,1.8){\line(1,-1){1}}
\put(6,.8){\line(1,1){.4}}
\put(7,1.8){\line(-1,-1){.4}}
\put(3.67,1.3){$\longrightarrow$}
\put(3.25,0){\small Figure 8}
\end{picture}}

\smallskip

In the double cover, $Y_B$, the solid torus $S^1\x D^2$ lifts to a bundle
$V$ over the circle with fiber the double cover of $D^2$ branched over
$2m$ points, a twice-punctured surface $S$ of genus $m-1$. The path
$\z$ lifts to a simple closed loop $\tz\C S\C  Y_B$, and changing $B$
by a  half-twist of along $\z$ as described corresponds to changing the
monodromy of the lifted bundle by a single Dehn twist along $\tz$. (This
is true essentially because each half-twist along $\z$ lifts to a full
twist in the double cover.)

Thus if $B'$ is the braid with the new positive half-twist, then its
corresponding double branched cover, $Y_{B'}$, is obtained from $Y_B$ by
cutting out $V$ and replacing it with the bundle over $S^1$ with fiber $S$
but whose monodromy is the monodromy of $V$ composed with a Dehn twist
about $\tz$. This means that $Y_{B'}$ is obtained by ($+1$)-Dehn surgery
on $\tz$ with respect to the  0--framing given by the pushoff of $\tz$
in the fiber $S$ of $V$. (For example, see \cite{ADK}.)

\begin{prop}\label{010}
The result of 0--surgery on $\tz$ in $Y_B$ is the double cover of $M^3$
branched along the link obtained from $B$ by the operation of Figure~9.
\end{prop}

\centerline{\unitlength 1cm
\begin{picture}(8,2.7)
\put(0,.8){\line(1,0){3}}
\put(0,1.8){\line(1,0){3}}
\put(1.5,.8){\line(0,1){1}}
\put(1.2,1.2){\small{$\zeta$}}
\put(5,.8){\line(1,0){1}}
\put(5,1.8){\line(1,0){1}}
\put(7,.8){\line(1,0){1}}
\put(7,1.8){\line(1,0){1}}
\put(6,1.8){\line(0,-1){1}}
\put(7,1.8){\line(0,-1){1}}
\put(3.67,1.3){$\longrightarrow$}
\put(3.25,0){\small Figure 9}
\end{picture}}

\begin{proof}
If we restrict the deck transformation $\t\co Y_B\to Y_B$ of the branched
cover to an annular neighborhood of $\tz$ in a fiber $S$ of $V$, then we
see an annulus double covering a disk with two branch points. Identify
a neighborhood of $\tz$ in $Y_B$ with $\tz\x I\x I$ where $I=[-1,1]$.
The restriction of $\t$ to this neighborhood is equivalent to
$\t(z,s,t)=(\bar{z},-s,t)$, and its fixed point set consists of two arcs
$\{(\pm 1, 0)\}\x I$ (identifying $\tz$ with $S^1$).
If we now change coordinates so that $I\x I$ becomes $D^2\C\bf{C}$,
then we get $\t(z,w)=(\bar{z},\rho{w})$, where $\rho$ is reflection in
the imaginary axis, and the fixed set is $\{\pm1\}\ \x$ \{the imaginary
axis $\cap\ D^2$\}.

According to our framing convention, 0--surgery on $\tz$ is the one
that kills the homology class of a pushoff of $\tz$ in $S$, ie the
class of $\tz\x \pt$ in $\tz\x D^2$. Thus the result of 0--surgery is
$Z=Y_B\- (\tz\x D^2)\cup_\vp(S^1\x D^2)$ where $\vp_*[\pt\x\bd
D^2]=[\tz\x \pt]$. Such a map $\vp$ is given by $\vp(z,w)=(w,z)$. Define
the involution $\s$ on $S^1\x D^2$ by $\s(z,w)=(\rho(z),\bar{w})$. Then
we see that the diagram
\begin{equation}\label{cd2}
   \begin{CD} 
 {S^1\x\bd D^2} @> {\s} >> { S^1\x\bd D^2} \\ 
  @ V{\vp} VV  @ VV {\vp} V \\ 
 {\tz\x \bd D^2} @>> {\t} > {\tz\x \bd D^2}
 \end{CD} 
\end{equation}
commutes. Thus, the restriction of $\t$ to $Y_B\- (\tz\x D^2)$ extends
to an involution $\t'$ over all of $Z$ via $\s$.

On the solid torus $S^1\x D^2$ the fixed set of $\t'=\s$ is
$\{\pm i\}\ \x \{\text{the real axis}~\cap\ D^2$\}. Thus the picture in
the quotient is exactly that of Figure~9.
\end{proof}

We now apply this proposition to the case at hand, where the 3--manifold
is $M_K$, the braid is the trivial braid with components $m_1$ and
$m_2$, and the arc $\z$ is the path $P$ of Figure~5. It follows that
$X_T(0,1,0)$ is the double branched cover of $E(1)_K$ with branch set
$S^1\x C$ where $C$ is the loop shown in Figure~10. Notice that $C$
is an unknotted circle which is unlinked from $K$.

\centerline{\unitlength 1cm
\begin{picture}(5,5.75)
\put (3,3){\oval(3,4)[r]}
\put (2,5){\line(3,-4){.4}}
\put (3,3.66){\line(-3,-4){1}}
\put (3,3.66){\line(-3,4){.4}}
\put (2,3.66){\line(3,4){1}}
\put (2,3.66){\line(3,-4){.4}}
\put (3,2.33){\line(-3,4){.4}}
\put (2,2.33){\line(3,-4){.4}}
\put (3,2.33){\line(-3,-4){1}}
\put (3,1){\line(-3,4){.4}}
\curve(2,5,1,5,0.6,4.75,0.5,4)
\curve(.47,2.25,.45,2.25)
\curve(.4,4.3, .35,4.29, .3,4.27, .25,4.22, .225,4.15, .224,4.05, .275,3.95, .35,3.92, .4,3.915) 
\curve(.58,3.915, .645, 3.95, .696,4.05, .695,4.15, .67,4.22,.62,4.27, .57,4.3 )
\put( .58,3.915){\line(0,-1){2.05}}
\put(.4,3.915){\line(0,-1){2.05}}
\curve( .35,2.24, .3,2.22, .25,2.17, .225,2.1, .224,2, .275,1.9, .35,1.87, .4,1.865) 
\curve(.58,1.865, .645, 1.9, .696,2, .695,2.1, .67,2.17,.62,2.22)
\curve(2,1,1,1,0.6,1.25,.5,1.75)
\curve(.5,1.6,.5,4)
\curve(.55,2.25,.57,2.24)
\put(4.25,4.9){\small $0$}
\put(.7,3.8){\small $C$}
\put (1.85,.5){\small Figure 10}
\end{picture}}

The double cover of a 3--ball, branched over an unknot, is $S^2\x I$
so it follows that the double cover of $M_K$ branched over $C$ is $M_K\#
M_K$.  Thus 
$$X_T(0,1,0) = E(1)\#_F \left(S^1\x (M_K\# M_K)\right)\#_F E(1)$$
This means that $X_T(0,1,0)$ is split by $S^1\x S^2$ with $b_2^+$
positive on each side. It follows that $\sw_{X_T(0,1,0)}=0$.

Next we need to make a similar calculation for $X_{T'}(0,1,0)$. This
time the arc $\z$ is the path $P'$ of Figure~5, and $X_{T'}(0,1,0)$
is the double branched cover of $E(1)_K$ with branch set $S^1\x C'$
where $C'$ is the loop shown in Figure~11.

\centerline{\unitlength 1.25cm
\begin{picture}(5.5,5.25)
\put (3,3){\oval(3,4)[r]}
\put (2,5){\line(3,-4){.4}}
\put (3,3.66){\line(-3,-4){1}}
\put (3,3.66){\line(-3,4){.4}}
\put (2,3.66){\line(3,4){1}}
\put (2,3.66){\line(3,-4){.4}}
\put (3,2.33){\line(-3,4){.4}}
\put (2,2.33){\line(3,-4){.4}}
\put (3,2.33){\line(-3,-4){1}}
\put (3,1){\line(-3,4){.4}}
\curve(2,5,1,5,0.6,4.75,0.5,4)
\curve(.4,4.3, .35,4.29, .3,4.27, .25,4.22, .225,4.15, .224,4.05, .275,3.95, .35,3.92, .4,3.915, .52,3.915, .57,3.92, .645, 3.95,.7,4.05)
\curve(.4,2.25, .35,2.24, .3,2.22, .25,2.17, .225,2.1, .224,2,.275,1.9, .35,1.87, .4,1.865, .52,1.865, .57,1.87, .645, 1.9, .696,2)
\curve(.65,2.18,.55,2.23)
\curve(2,1,1,1,0.6,1.25,.5,1.75)
\curve(.4,3.27, .35,3.26, .3,3.24, .25,3.19, .225,3.12, .224,3.02, .275,2.92, .35,2.89, .4,2.885, .52,2.885, .57,2.89, .645, 2.92, .696,3.02, .695,3.12, .67,3.19,.62,3.24, .57,3.27 )
\put(0,2.72){\small{$m_0$}}
\curve(.5,3.8,.5,3.02)
\curve(.5,2.8,.5,2)
\curve(.645,4.25,.58,4.27)
\curve(.645,4.25,1.5,4.4,2.4, 4.37)
\curve(.7,4.05,1.5,4.2,2.4, 4.25)
\curve(2.62, 4.25,3.4, 3.67,2.55, 3.02)
\curve(2.55,4.37,2.8,4.3, 3, 4.16)
\curve(3.04, 4.03,3.1,3.5,2.6,3.1)
\curve(2.4, 2.94,2.2,2.85)
\curve(2.4,3.05,2.2,2.95,2,2.7,1.8,2.6)
\curve(2.1,2.92,1.95,2.9,1.8,2.7,1.6,2.45,.7,2)
\curve(1.7,2.65,1,2.35, .65,2.18)
\put(4.05,4.6){\small$0$}
\put(1,1.95){\small$C'$}
\put (2,.5){\small Figure 11}
\end{picture}}

\centerline{\unitlength 1.5cm
\begin{picture}(9,4)
\put(1,4){\line(1,0){3}}
\put(1.5,3.5){\line(1,0){2}}
\put(1,4){\line(0,-1){1.5}}
\put(1.5,3.5){\line(0,-1){1}}
\curve(1.5,2.45,1.7,2.6,1.55,2.8)
\curve(1,2.45,.8,2.6,.95,2.8)
\curve(1,2.45,1.25,2.4,1.5,2.45)
\curve(1.07,2.85,1.25,2.87,1.43,2.85)
\put(1,2.4){\line(0,-1){.85}}
\put(1,.9){\line(0,1){.5}}
\put(1.5,2.4){\line(0,-1){1}}
\put(1,.9){\line(1,0){1.3}}
\put(1.5,1.4){\line(1,0){.8}}
\put(2.3,.9){\line(1,1){.5}}
\curve(2.3,1.4,2.5,1.2)
\curve(2.8,.9,2.6,1.1)
\put(2.8,1.4){\line(1,0){.7}}
\put(2.8,.9){\line(1,0){1.2}}
\put(3.5,1.4){\line(0,1){.8}}
\put(4,.9){\line(0,1){1.3}}
\curve(3.5,2.2,3.53,2.33,3.57,2.45)
\curve(3.65,2.52,3.75,2.55,3.9,2.45,4,2.2)
\put(3.5,3.5){\line(0,-1){.8}}
\put(4,4){\line(0,-1){1.4}}
\curve(3.5,2.7,3.53,2.6,3.6,2.5,3.75,2.4,3.86,2.41)
\curve(4,2.6,3.97,2.5,3.93,2.45)
\curve(1.05,1.78,1.2,1.6,1,1.45,.8,1.6,.95,1.78)
\put(.5,1.7){\small$m_0$}
\put(.5,2.6){\small$C'$}
\put(4.2,3.7){\small$0$}
\put(2,.3){\small Figure 12}
\put(6,3.5){\line(1,0){.8}}
\put(5.5,4){\line(1,0){1.3}}
\put(6.8,3.5){\line(1,1){.5}}
\curve(6.8,4,7,3.8)
\curve(7.3,3.5,7.1,3.7)
\put(7.3,3.5){\line(1,0){.7}}
\put(7.3,4){\line(1,0){1.2}}
\put(8,3.5){\line(0,-1){.8}}
\put(8.5,4){\line(0,-1){.75}}
\put(8.5,3.15){\line(0,-1){.55}}
\curve(8.55,3.53,8.7,3.35,8.5,3.2,8.3,3.35,8.45,3.53)
\put(8.7,3.45){\small$m''_0$}
\curve(8,2.7,8.03,2.6,8.1,2.5,8.25,2.4,8.36,2.41)
\curve(8.5,2.6,8.47,2.5,8.43,2.45)
\put(6,3.5){\line(0,-1){.9}}
\put(5.5,4){\line(0,-1){1.3}}
\curve(5.5,2.7,5.53,2.6,5.6,2.5,5.75,2.4,5.86,2.41)
\curve(6,2.6,5.97,2.5,5.93,2.45)
\curve(5.5,2.2,5.53,2.33,5.57,2.45)
\curve(5.65,2.52,5.75,2.55,5.9,2.45,6,2.2)
\curve(8,2.2,8.03,2.33,8.07,2.45)
\curve(8.15,2.52,8.25,2.55,8.4,2.45,8.5,2.2)
\put(5.5,.9){\line(0,1){.5}}
\put(5.5,1.5){\line(0,1){.7}}
\curve(5.55,1.78,5.7,1.6,5.5,1.45,5.3,1.6,5.45,1.78)
\put(4.98,1.68){\small$m'_0$}
\put(6,1.4){\line(0,1){.8}}
\put(8,1.4){\line(0,1){.8}}
\put(8.5,.9){\line(0,1){1.3}}
\put(6,1.4){\line(1,0){.8}}
\put(5.5,.9){\line(1,0){1.3}}
\put(6.8,.9){\line(1,1){.5}}
\curve(6.8,1.4,7,1.2)
\curve(7.3,.9,7.1,1.1)
\put(7.3,1.4){\line(1,0){.7}}
\put(7.3,.9){\line(1,0){1.2}}
\put(6.1,2.9){\small$-2$}
\put(7.6,1.9){\small$-2$}
\put(6.5,.3){\small Figure 13}
\end{picture}}

An isotopy of Figure~11 gives Figure~12, and the corresponding
double branched cover is shown in Figure~13.  (For techniques for
determining this double cover, see \cite{Rolfsen}.) Using Fox calculus,
one calculates the torsion of the link $m'_0\cup m''_0$ in the double
cover $(M_K)_{C'}$ (the 3--manifold shown in Figure~13). This torsion is
$1$. According to \cite{MT} this means that the Seiberg--Witten invariant
of $S^1\x \left((M_K)_{C'}\- (m'_0\cup m''_0)\right)$ is $1$. Since the
Seiberg--Witten invariant of $E(1)\- F$ is also equal to $1$, the gluing
theorem of Taubes \cite{T} tells us that the Seiberg--Witten invariant of
$E(1)\#_{F=S^1\x m'_0} (S^1\x {(M_K)}_{C'})\#_{S^1\x m'_0=F} E(1)$ 
is equal to $1$. However, we have just shown that this manifold is
$X_{T'}(0,1,0)$.

\begin{prop}
For the nullhomologous Lagrangian tori $T$, $T'$ in $X$, we have
$$\sw_{X_T(0,1,0)}=0, \quad \sw_{X_{T'}(0,1,0)}=1$$
\end{prop}

\section{Calculation of $I(X,T)$: $X(1,0,0)$}\label{L}

The key calculation of this section will show that the Seiberg--Witten
invariants of the manifolds $X_T(1,0,0)$ and $X_{T'}(1,0,0)$ vanish. Our
approach here is to describe the surgered manifolds in terms of a
branched covering. (It would be useful to compare with \cite{Lagr},
where a more general approach is utilized.)

\begin{prop}\label{100}
Let $\bg$ denote either $\g$ or $\g'$, and let $Z$ be the result of the
surgery on $T=S^1\x \bg\C X$ which kills $S^1\x\pt \x \pt \C S^1\x\bg\x\bd
D^2\C S^1\x \left(M_{K\# K}\- (m'_0\cup m''_0)\right)$. Then $Z$ is
the double branched cover of the manifold $W$ obtained from $E(1)_K$
by a surgery on a circle $S^1\x \{\text{point on $\bg$}\} $ (trading a
neighborhood $S^1\x D^3$ for $D^2\x S^2$). The branch set in $W$ of this
cover consists of a pair of disjoint 2--spheres of self-intersection 0.
\end{prop}
\begin{proof}
As we have seen in the proof of Proposition~\ref{010}, the deck
transformation of $X\to E(1)_K$ in a neighborhood $S^1\x\bg\x D^2$ of $T$
is given by $\t(t,z,w)=(t,\bar{z},\rho(w))$ where $\rho$ is reflection
through the imaginary axis. The manifold $Z$ is:
$$Z=\left(X\- (S^1\x\bg\x D^2)\right)\cup_{\vt}(S^1\x S^1\x D^2)$$
where $\vt(t,z,w)=(w,z,t)$. Then the diagram
\begin{equation}\label{cd3}
   \begin{CD} 
 {S^1\x S^1\x\bd D^2} @> {\upsilon} >> {S^1\x S^1\x\bd D^2} \\ 
  @ V{\vt} VV  @ VV {\vt} V \\ 
 {S^1\x\bg\x \bd D^2} @>> {\t} > {S^1\x\bg\x \bd D^2}
 \end{CD} 
\end{equation}
commutes, where $\upsilon(t,z,w)=(\rho(t),\bar{z},w)$. Thus $\upsilon$
extends the deck transformation $\t$ over the surgered manifold $Z$.

The quotient $(S^1\x S^1\x D^2)/\t\cong S^1\x(D^2\x I)\cong S^1\x
D^3$, but
$(S^1\x S^1\x D^2)/\upsilon\cong S^2\x D^2$, since the action of
$\upsilon$ restricted to $S^1\x S^1\x\{\pt\}$ is equivalent to the
action of the deck transformation of the double covering $T^2\to S^2$
with four branch points.
Thus the effect of the surgery on the base is to perform surgery on the
circle $S^1\x\{\pt\}\C S^1\x D^3$. Before performing the surgery, the
branch set consists of two tori. Since the fixed point set of $\upsilon$
on $S^1\x S^1\x D^2$ is $\{\pm i\}\x\{\pm1\}\x D^2$, the surgery trades
a pair of annuli for four disks. Removing the annuli leaves us with a
pair of complementary annuli in the branch set, and the addition of the
four disks caps them off, giving a pair of 2--spheres.

To see that the components of the branch set of $W$ have self-intersection
$0$, first consider the branch torus $S^1\x m_1$ of $X$. Write
$m_1=J_1\cup J_2$, the union of two intervals meeting only at their
endpoints. We do this so that the intersection of $S^1\x m_1$ with $(S^1\x
S^1\x D^2)/\t\cong S^1\x D^3$ is $S^1\x J_2$ where $J_2\cap \bd D^3 =
\bd J_2$. Then the corresponding component of the branch set in $W$
is $(S^1\x J_1) \cup (D^2\x \bd J_2)$. In $X$ we can isotop $S^1\x m_1$
slightly by moving $m_1$ to $\bar{m}_1 = I_1\cup I_2$ in $M_K$; so that
$\bar{m_1}\cap m_1 = \emptyset$, $\bar{m}_1\cap D^3 = I_2$, and $I_2
\cap \bd D^3 = \bd I_2$. Then in $W$, $(S^1\x I_1)\cup (D^2\x \bd I_2)$
is disjoint from $(S^1\x J_1)\cup (D^2\x \bd J_2)$.
\end{proof}

Let $\G_i$ denote the components of the branch set in $Z$. The $\G_i$
are also 2--spheres of self-intersection 0.

In $X$ there is a `section class' $C$ which arises from the sections of
the elliptic fibrations on the copies of $E(1)$. To build a representative
for $C$, start with a fixed Seifert surface $B_0$ of $K\# K$ whch is a
fiber of the fibration of $S^3\- K\# K \to S^1$. The boundary of $B_0$
is capped off by the 2--disk introduced when we do 0--surgery on $K\#
K$ to form $M_{K\# K}$. The tori $S^1\x m'_0$ and $S^1\x m''_0$ which are
identified with fibers of $E(1)'$ and $E(1)''$  each intersect $\{\rm pt\}
\x B_0$ in a single point. Remove disks in $\{\rm pt\} \x B_0$ about each
of these points. The boundaries then bound disks of self-intersection
$-1$, sections of $E(1)$ minus the neighborhood of a fiber. The union
of these surfaces gives a genus--2 surface of self-intersection $-2$
representing $C$.

The loop $\bg\C X$ is contained in a Seifert surface for $K\# K$, and
we may assume that it is disjoint from $B_0$. Thus the surgery torus,
$S^1\x \bg$ is disjoint from $C$.
Since $C\cdot (S^1\x m_1) = C\cdot F' = 1$, after surgery in Z, we
still have $C\cdot \G_1=1$; so $\G_1$ is an essential 2--sphere in $Z$
(and similarly for $\G_2$).  Thus $Z$ contains an essential 2--sphere
of self-intersection 0, and that means that  $\sw_Z=0$ \cite{Turkey}.

\begin{thm}\label{T}
$T$ and $T'$ are essential and homologous Lagrangian tori of $X$;
however, there is no orientation-preserving diffeomorphism $f$ of $X$
with $f(T)=T'$.
\end{thm}
\begin{proof}
We have $\sw_{X_T(1,0,0)}=0$, $\sw_{X_T(0,1,0)}=0$, and, since
$X_T(0,0,1)=X$,  $\sw_{X_T(0,0,1)}=(t_F^2-1+t_F^{-2})^2$. Hence 
$I(X,T)=\{r(t_F^2-1+t_F^{-2})^2|\,  r\in \Z\}$. On the other hand,
$\sw_{X_{T'}(1,0,0)}=0$, $\sw_{X_{T'}(0,1,0)}=1$, and
$X_{T'}(0,0,1)=X$; so $I(X,T')=\{q+r(t_F^2-1+t_F^{-2})^2| \, q,r\in \Z\}$.
This concludes the proof since $I(X,T)$ is an orientation-preserving
diffeomorphism invariant of $(X,T)$.
\end{proof}

Auroux, Donaldson, and Katzarkov have shown in \cite{ADK}
that  the surgery manifolds $X_T(0,k,1)$ and $X_{T'}(0,k,1)$ are
symplectic for all $k\in\Z$. The corresponding Seiberg--Witten
invariants are $\sw_{X_T(0,k,1)}=(t_F^2-1+t_F^{-2})^2$ and
$\sw_{X_{T'}(0,k,1)}=k+(t_F^2-1+t_F^{-2})^2$. Note that the leading
coefficient of these polynomials is $\pm 1$, as required by Taubes'
theorem.

\Addresses
\end{document}